\documentclass[a4paper,10pt]{article}

\usepackage[a4paper,left=3cm,right=3cm,top=2.5cm,bottom=2.5cm]{geometry}
\usepackage{hyperref}
\usepackage{subfig}
\usepackage{pgfplots}
\usepackage{pgfplotstable}
\usepackage{mathtools}
\usepackage{multicol}
\usepackage{comment}
\usepackage{booktabs}
\pgfplotsset{compat=1.5}
\usepackage{amssymb}
\usepackage{url}
\usepackage{bm}

\usepackage{pdfsync}
\usepackage{float}
\usepackage{tabularx}
\usepackage{enumerate}
\usepackage{array}
\usepackage{xspace}
\usepackage{siunitx}
\usepackage{authblk}

\usepackage{amsthm}
\usepackage{amsmath}
\usepackage{stmaryrd}
\usepackage[ruled,vlined,linesnumbered]{algorithm2e}
\usepackage{graphicx}
\usepackage{booktabs}
\usepackage{cleveref}

\usepackage[normalem]{ulem}

\newtheorem{example}{Example}

\DeclareMathOperator*{\argmin}{\arg\!\min}

\begin{document}

\title{Data-driven Discovery of Delay Differential Equations\\ with Discrete Delays}

\author[2]{Alessandro Pecile\footnote{pecile.alessandro@spes.uniud.it}}
\author[1]{Nicola Demo\footnote{ndemo@sissa.it}}
\author[3]{Marco Tezzele\footnote{marco.tezzele@austin.utexas.edu}}
\author[1]{Gianluigi Rozza\footnote{grozza@sissa.it}}
\author[2]{Dimitri Breda\footnote{dimitri.breda@uniud.it}}

\affil[1]{Mathematics Area, mathLab, SISSA, via Bonomea 265, I-34136 Trieste,
  Italy}
\affil[2]{CDLab – Computational Dynamics Laboratory 
Department of Mathematics, Computer Science and Physics – University of Udine,
via delle Scienze 206, 33100 Udine, Italy}
\affil[3]{Oden Institute for Computational Engineering and Sciences, University of Texas
at Austin, Austin, 78712, TX, United States}

\maketitle

\begin{abstract}
The Sparse Identification of Nonlinear Dynamics (SINDy) framework is a robust method for identifying governing equations, successfully applied to ordinary, partial, and stochastic differential equations. In this work we extend SINDy to identify delay differential equations  by using an augmented library that includes delayed samples and Bayesian optimization. To identify a possibly unknown delay we minimize the reconstruction error over a set of candidates. The resulting methodology improves the overall performance by remarkably reducing the number of calls to SINDy with respect to a brute force approach. We also address a multivariate setting to identify multiple unknown delays and (non-multiplicative) parameters. Several numerical tests on delay differential equations with different long-term behavior, number of variables, delays, and parameters support the use of Bayesian optimization highlighting both the efficacy of the proposed methodology and its computational advantages. As a consequence, the class of discoverable models is significantly expanded.
\end{abstract}

\tableofcontents

\section{Introduction}
\label{sec:intro}
In the last years data-driven model discovery has enhanced our ability to understand, forecast, and control dynamical systems by combining techniques from the abstract theory of dynamical systems, data science, and numerical analysis. Even if a lot of research has been devoted to the identification of Ordinary Differential Equations (ODEs) or Partial Differential Equations (PDEs), little to none attention has been reserved to functional equations, among which stands the class of Delay Differential Equations (DDEs). In this work we leverage data-driven system identification methods and Bayesian optimization (BO) to efficiently compute the unknown delays and possibly other parameters of an unknown governing DDE.

Different techniques have been proposed to approximate or even reconstruct a dynamical system from data. To name a few, dynamic mode decomposition \cite{tu2014dmd, kutz2017,schmidt2022}, PDEs discovery through sparse optimization~\cite{schaeffer2017learning}, sparse identification of nonlinear dynamics (SINDy)~\cite{Brunton2016, Brunton2019, Champion2019}, physics-informed neural networks \cite{Raissi2019, Karniadakis2021, Cuomo2022, Sholokhov2023}, and operator inference~\cite{peherstorfer2016data, qian2020lift, mcquarrie2021data} represent today consolidated approaches. In particular, SINDy was originally proposed in \cite{Brunton2016} to reconstruct an autonomous ODE $x' = f(x)$ from a time series of the state variable $x$. The right-hand side $f$ is recovered by assuming it lies in the linear span of a chosen library of functions, reducing the problem to a linear regression that is solved by imposing sparsity constraints (owing to the observation that typically $f$ is made of few terms). The use of SINDy has spread quite rapidly and now several tools and packages are available which extend to, e.g., partial \cite{Rudy2017, schaeffer2017learning, Kaheman2020} or stochastic ordinary \cite{Boninsegna2018, Huang2022, WANG2022, jacobs2023hypersindy} differential equations (see also \cite{wanner2024higher} and the references therein). In particular, in this work we resort to the open source Python package \texttt{pySINDy}~\cite{deSilva2020, Delahunt2021,Kaptanoglu2022}.

Although SINDy and its variations have proven effective in different settings, from ODEs \cite{Brunton2016, Kaheman2020} to plasma physics \cite{Kaptanoglu2021mhd, Yifei2021, Alves2022} and fluid dynamics \cite{LoiseauNoackBrunton2018, LoiseauBrunton2018, Loiseau2020, Kaptanoglu2021gs, Fukami2021, Oishi2023}, several models used in applications are not yet fully covered due to their specific mathematical formulation, such as DDEs. DDEs express the current-time derivative of the unknown function through its history. In their simplest form, such history spans a finite time interval and is given by values on fixed time instants in the past, known as the case of (constant) {\it discrete} delays (opposed to the case of {\it distributed} delays, where the past dependence is represented by an integral over some past interval). DDEs are useful to model phenomena characterized by a reaction time or latency in a response, justifying their massive use in control theory (to have a realistic model of the action of the controller) and in epidemiology (to effectively model the time between the exposure to a disease, the manifestation of the symptoms and healing). The mathematical treatment of DDEs \cite{Hale1977, Diekmann1995, erneux2009} requires tools from functional analysis, such as semigroup and generation theory, that make it possible to express their properties. Among these properties emerges the intrinsic infinite dimension of their {\it dynamical} state (the object evolving in time, opposed to the finite dimension of their {\it physical} state, i.e., the number of unknown variables).

To the best of our knowledge preliminary results using SINDy were presented by the authors at \cite{breda2022data}, while a first publication~\cite{sandoz2023sindy} appeared as recently. The two proposed approaches are similar, both based on the idea of extending the library of functions candidate to form the right-hand side by including delayed terms besides the current-time samples. Both are limited to the case of a single constant discrete delay and result effective in identifying the latter when unknown, by selecting among a set of candidates the value that minimizes the reconstruction error returned by SINDy. Additionally,~\cite{sandoz2023sindy} successfully applies the technique to realistic data to recover a plausible model for bacterial zinc response, thus highlighting the potentiality of the proposed approach.

\begin{figure}[ht]
\centering
\includegraphics[width=\textwidth]{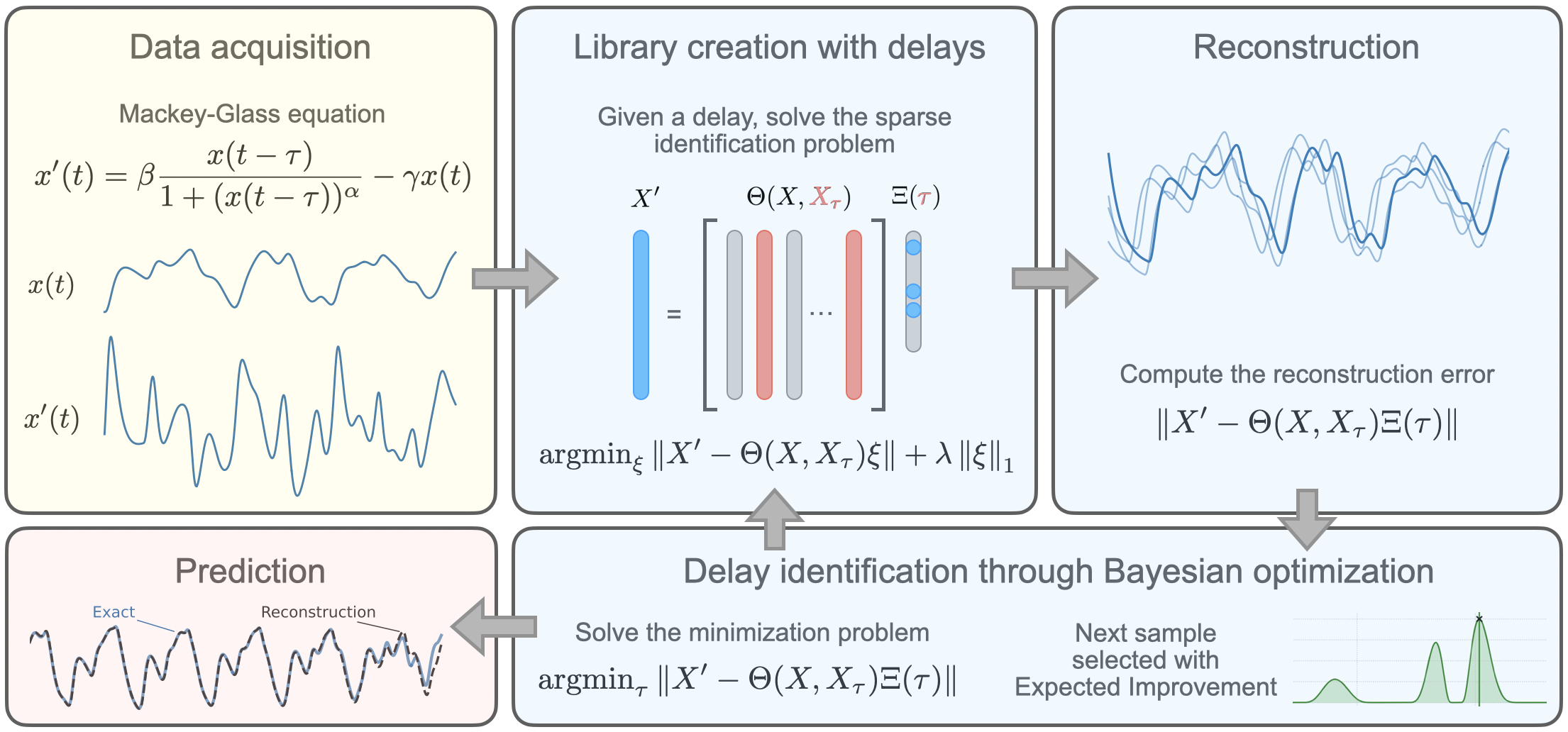}
\caption{Abstract representation of the information flow to identify DDEs with SINDy. First we collect the time series generated from an unknown process (yellow box), then we enter in the optimization loop (blue boxes). The optimization starts with a candidate delay $\tau$, which is used to form the functions library for the sparse identification problem. The BO loop minimizes the reconstruction error until convergence. The optimized model correctly identifies the unknown delay and can be used for future state predictions (red box).}\label{fig:scheme}
\end{figure}

In this paper our aim is twofold. On the one hand, we propose to use BO to identify the unknown delay, instead of the brute force approach of spanning a whole set of candidate values. This largely improves the performance of the overall methodology in terms of the number of SINDy reconstructions. On the other hand, a multivariate BO approach is proposed to identify multiple unknown (non-multiplicative) parameters. Consequently, the class of discoverable models is considerably enlarged, including, e.g., systems with specific nonlinear terms like peculiar functional responses in biological models, as well as systems with multiple constant discrete delays. A sketch of the end-to-end pipeline is depicted in Figure~\ref{fig:scheme}, where we consider the Mackey-Glass equation, one of the test cases analyzed in Section~\ref{sec:results}.

The structure of the paper is the following. In Section~\ref{sec:dde} we briefly introduce DDEs, together with four test cases. In Section~\ref{sec:sindy} we summarize the SINDy algorithm for ODEs and its extension to DDEs with delays assumed to be known. In Section~\ref{sec:bopt} we recall the basics of BO in order to describe in Section~\ref{sec:methods} our general methodology for applying SINDy to DDEs with multiple unknown delays and/or parameters. In Section~\ref{sec:results} we describe and discuss the numerical results obtained for the four test cases. In Section~\ref{sec:conclusion} we draw some conclusions and propose future developments.

\section{Delay Differential Equations}
\label{sec:dde}
To introduce DDEs we follow~\cite{Breda2015}. Here we summarize the basic elements and the relevant notation.

Let $\tau>0$ be real, $n$ be a positive integer and $C\coloneqq C([-\tau,0],\mathbb{R}^{n})$ be the Banach space of continuous functions $[-\tau,0]\to\mathbb{R}^{n}$ endowed with the uniform norm
\begin{equation*}
\left\|\varphi\right\|_{C}\coloneqq\max_{s\in[-\tau,0]}\left\|\varphi(s)\right\| ,
\end{equation*}
for $\|\cdot\|$ any norm in $\mathbb{R}^{n}$. A generic (nonlinear, autonomous) DDE is an equation of the form
\begin{equation}
\label{eq:ddeDB}
x'(t)=f(x_{t}),\quad t\geq0,
\end{equation}
where $f:C\to\mathbb{R}^{n}$, $x_{t}\in C$ is defined as
\begin{equation}\label{eq:state}
x_{t}(s)\coloneqq x(t+s),\quad s\in[-\tau,0],
\end{equation}
and $\tau$ is the (maximum) delay. An Initial Value Problem (IVP) for \eqref{eq:ddeDB} is obtained by assigning an initial history through a given function $\varphi\in C$, i.e., by imposing
\begin{equation*}
x(s)=\varphi(s),\quad s\in[-\tau,0].
\end{equation*}
Under standard smoothness assumptions on $f$ (e.g., Lipschitz continuity, see \cite{Hale1977}), the IVP has a unique solution which depends smoothly on $\varphi$, thus giving rise to a dynamical system on $C$ where the state at time $t\geq0$ is represented by the relevant history $x_{t}$ as defined in \eqref{eq:state}. Consequently, DDEs originate infinite-dimensional dynamical systems despite their finite physical dimension $n$. Precisely this infinite dimension provides room for complex dynamical behaviors already in the scalar case \cite{Sprott2007}.

Although the setting introduced above is favorable in view of generality, the cases of interest in this work belong to the class of DDEs with a finite number $k$ of constant discrete delays $0<\tau_{1}<\cdots<\tau_{k}$, i.e., equations of the form
\begin{equation}\label{eq:ddediscrete}
x'(t)=f(x(t),x(t-\tau_{1}),\ldots,x(t-\tau_{k})),
\end{equation}
where now, with slight abuse of notation, $f:\mathbb{R}^{n(k+1)}\to\mathbb{R}^{n}$. A brief description of the numerical test cases considered in Section \ref{sec:results} follows and completes this section.

\begin{example}[Delay logistic equation]\label{ex:logistic}
The delay logistic or Hutchinson equation \cite{Hutchinson1948} is among the simplest DDEs. It reads
\begin{equation}\label{eq:logistic}
x'(t)=\rho x(t)(1-x(t-\tau)),
\end{equation}
with $\rho>0$.
The behavior of the solution to the associated IVP can be described in terms of $\rho$ and $\tau$: monotonic for $\rho\tau<1/e$, oscillatory for $1/e<\rho\tau<\pi/2$. In both cases, the solution approaches the equilibrium $\bar x=1$ as $t\to+\infty$. At $\rho\tau=\pi/2$ a Hopf bifurcation takes place and for $\rho\tau>\pi/2$ the solution approaches a limit cycle. It is possible to express the limit cycle using the Lambert $W$ function \cite{Gopalsamy1992}. The profile of the cycle changes with increasing $\rho$ and leads to regions with large derivatives separated by plateaus~\cite{Gopalsamy1992}.
\end{example}

\begin{example}[SIR model with delay]\label{ex:sir}
To test the proposed approach on a system ($n>1$) of DDEs we use a modified version of the celebrated SIR model used in mathematical epidemiology to study the diffusion of a disease and first appeared in \cite{Kermack1927} as a trivial instance of a more complex integral equation with age structure. The modification concerns in introducing a delay in the recovered individuals leaving the infective compartment, thus modeling an average healing time:
\begin{equation}\label{eq:sir}
\left\{\setlength\arraycolsep{0.1em}\begin{array}{rcl}
s'(t)&=&-\beta s(t)i(t),\\[1mm]
i'(t)&=&\beta s(t)i(t)-\mu i(t-\tau),\\[1mm]
r'(t)&=&\mu i(t-\tau).
\end{array}\right.
\end{equation}
Above $\beta$ and $\mu$ are positive (multiplicative) parameters. This model, considered also in \cite{Ebraheem2021} in the context of optimal control, lacks the natural positivity of solutions required for a population model. Yet it expresses a range of behaviors wider than that of the classical SIR model and it is thus good to test the proposed identification strategy in the case of systems. For other SIR models with delay see \cite{breda2024viability} and the references therein.
\end{example}

\begin{example}[Mackey-Glass equation]\label{ex:mg}
The Mackey-Glass equation~\cite{Mackey1977} is the nonlinear DDE
\begin{equation}\label{eq:mg}
x'(t)=\beta\frac{x(t-\tau)}{1+(x(t-\tau))^{\alpha}}-\gamma x(t),\quad
\end{equation}
with $\beta$ and $\gamma$ positive reals and $\alpha$ a positive integer. It was first proposed to describe the production of blood cells in physiology, based on a feedback system relying upon a retarded Hill-type response with Hill coefficient $\alpha$, i.e., the term
\begin{equation}\label{eq:hill}
h(t-\tau;\alpha)\coloneqq\frac{1}{1+(x(t-\tau))^{\alpha}}.
\end{equation}
Depending on the values of the parameters, \eqref{eq:mg} displays a range of periodic and possibly chaotic dynamics with chaotic attractors. Indeed, it is a concrete example of how scalar DDEs can show dynamically rich behaviors that are excluded for scalar or planar ODEs. In Section~\ref{sec:results} we will apply our method to simultaneously identify both $\tau$ and $\alpha$.
\end{example}

\begin{example}[Multiple delays]\label{ex:multiple}
As a case with more than one delay we consider the DDE
\begin{equation}\label{eq:multiple}
x'(t)=ax^2(t-\tau_{1})+bx^3(t-\tau_{2}) ,
\end{equation}
from \cite{Ji2021learning}. Therein it is used to test the application of time delay neural networks (neural networks with delayed input \cite{zhu2023neural, zhu2021neural}) to learn DDEs while simultaneously training the unknown delays. For the use of neural networks to learn DDEs see also \cite{Ji2024Learn}.
\end{example}

\section{Sparse Identification of Nonlinear Dynamics}
\label{sec:sindy}
SINDy~\cite{Brunton2016, Brunton2019, Champion2019} is a technique to discover dynamical systems models from data and observations. In its standard setting for ODEs, one supposes to have measurements of the unknown $x(t)$ of an ODE $x'=f(x)$ at time instants $t_{1},\ldots,t_{m}$ and assumes that each component $f_i$, $i=1\ldots,n$, of $f$ can be (approximately) represented as a sparse linear combination of selected basis functions $\theta_1,\ldots,\theta_{p}$, namely
\begin{equation}\label{eq:fi}
f_i\approx\sum_{j=1}^{p}\theta_{j}\xi_{j,i}.
\end{equation}
To describe the way the right-hand side $f$ is reconstructed from the data, i.e., to recover the coefficients $\xi_{j,i}$, $j=1,\ldots,p$ in \eqref{eq:fi} for every $i=1,\ldots,n$, it is convenient to introduce the following compact notation.

We collect in a matrix $X\in\mathbb{R}^{m\times n}$ all the available measurements, with a row per time instant and a column per variable:
\begin{equation*}
    X \coloneqq \begin{bmatrix}
        x_1(t_1) & x_2(t_1) & \dots & x_n(t_1) \\
        x_1(t_2) & x_2(t_2) & \dots & x_n(t_2) \\
        \vdots & \vdots & & \vdots \\ x_1(t_m) & x_2(t_m) & \dots & x_n(t_m)
    \end{bmatrix}.
\end{equation*}
The entries of a matrix of derivatives $X'\in\mathbb{R}^{m\times n}$ are organized similarly, namely
\begin{equation*}
    X' \coloneqq \begin{bmatrix} x_1'(t_1) & x_2'(t_1) & \dots & x_n'(t_1) \\
        x_1'(t_2) & x_2'(t_2) & \dots & x_n'(t_2) \\
        \vdots & \vdots & & \vdots \\
        x_1'(t_m) & x_2'(t_m) & \dots & x_n'(t_m)
    \end{bmatrix},
\end{equation*}
with derivative values known a priori from other measurements or conveniently approximated \cite{Chartrand2011}.

A matrix $\Theta(X)\in\mathbb{R}^{m\times p}$ is then introduced to collect all the $p$ basis functions (one per column) evaluated at all the $m$ time instants (one per row). An example may include polynomial terms up to degree $d$, as well as other nonlinear terms possibly suggested by some physics knowledge of the system at hands (e.g., trigonometric functions for periodic forcing terms):
\begin{equation*}
\Theta(X)\coloneqq\begin{bmatrix}
1 & X & X^2 & \cdots & X^d & \cdots & \sin(X) & \cdots
\end{bmatrix} .
\end{equation*}
Above, $X^{d}$ denotes a matrix with column vectors given by all possible time-series of $d$-th degree polynomials in the unknown $x$, e.g., for $d=2$,
\begin{equation*}
    X^2 \coloneqq \begin{bmatrix}
        x_1(t_{1})^2 & x_1(t_{1})x_2(t_{1}) & \cdots & x_1(t_{1})x_n(t_{1})& x_2(t_{1})^2 & x_2(t_{1})x_3(t_{1}) & \cdots & x_n^2(t_{1}) \\
x_1(t_{2})^2 & x_1(t_{2})x_2(t_{2}) & \cdots & x_1(t_{2})x_n(t_{2}) & x_2(t_{2})^2 & x_2(t_{2})x_3(t_{2}) & \cdots & x_n^2(t_{2}) \\
        \vdots & \vdots & & \vdots & \vdots & \vdots & & \vdots\\
        x_1(t_{m})^2 & x_1(t_{m})x_2(t_{m}) & \cdots &x_1(t_{m})x_n(t_{m}) & x_2(t_{m})^2 & x_2(t_{m})x_3(t_{m}) & \cdots & x_n^2(t_{m})
    \end{bmatrix}.
\end{equation*}
Similarly,
\begin{equation*}
\sin(X)\coloneqq \begin{bmatrix}
\sin(x_1(t_{1})) & \sin(x_2(t_{1})) & \cdots & \sin(x_n(t_{1})) \\
\sin(x_1(t_{2})) & \sin(x_2(t_{2})) & \cdots & \sin(x_n(t_{2})) \\
\vdots & \vdots & \ddots & \vdots \\
\sin(x_1(t_{m})) & \sin(x_2(t_{m})) & \cdots & \sin(x_n(t_{m})) \\
    \end{bmatrix} ,
\end{equation*}
and the same holds for any other possible nonlinear function.

Of course, the choice of the basis functions is crucial. According to \cite{Brunton2016, Brunton2019}, reasonable choices in the context of a dynamical system originated from ODEs are polynomial functions, but different contexts and models can lead to different choices. In any case, the number of columns $p$ of the library matrix $\Theta(X)$ depends on the set of candidate class(es) of functions. For instance, a polynomial library of degree $d$ for an $n$-dimensional ODE amounts to $p={n+d\choose d}$ basis functions.

Finally, by collecting all the unknown coefficients $\xi_{j,i}$ of \eqref{eq:fi} in a vector $\xi_{i}\in\mathbb{R}^{p}$ for every $i=1,\ldots,n$, the problem of determining the latter is compactly written as
\begin{equation}\label{eq:sindyknown}
X'=\Theta(X)\Xi ,
\end{equation}
where $\Xi\in\mathbb{R}^{p\times n}$ is the matrix
\begin{equation*}
\Xi\coloneqq\begin{bmatrix}
\xi_{1} & \xi_{2} & \cdots & \xi_{n}
\end{bmatrix}.
\end{equation*}
The latter is solved column-wise as the optimization problem
\begin{equation}\label{eq:optknown}
\xi_i = \argmin_{\xi \in \mathbb{R}^p}\left(\left\|X_{i}'-\Theta(X)\xi\right\|+\lambda\left\|\xi\right\|_1\right) ,
\end{equation}
for every $i=1,\ldots,n$ via convex $\ell_1$-regularized sparse regression. Above $X_{i}'$ is the $i$-th column of $X'$, $\lambda$ represents a sparsity-promoting parameter and different techniques can be adopted to achieve an effective solution such as the sequential thresholded least-squares (STLS) algorithm \cite{Brunton2016} or the least absolute shrinkage and selection operator (LASSO) regression \cite{tibshirani96regression}. All the computations presented in Section~\ref{sec:results} rely on the open source Python package \texttt{pySINDy}~\cite{deSilva2020, Delahunt2021,Kaptanoglu2022}. In particular, \eqref{eq:optknown} is solved cia STLS.

\subsection{SINDy for DDEs with known delays}
\label{sec:knowndelay}
Extending SINDy to DDEs poses additional challenges since the current-time evolution depends on the past, even if one restricts to the form \eqref{eq:ddediscrete} with a single delay. However, a helpful yet not necessarily realistic assumption is that of considering the delay as known. In this section we deal with this case. The extension to the case of multiple known delays is straightforward and we omit the details. Instead, the problem of identifying {\it unknown} delays is tackled in Section~\ref{sec:methods} via BO techniques as introduced in Section~\ref{sec:bopt}.

A possible solution to the case of a single constant discrete and known delay $\tau$ is to look at the delayed variables as at further state variables. This solution leads to consider also the delayed values of the measured data besides the current-time ones stored in the matrix $X$. Delayed values are readily available only if measurements are taken at uniformly distributed time instants and $\tau$ is an integer multiple of the relevant sampling time. If this is not the case, a simple linear interpolation can be invoked to reconstruct past values. So let $X_\tau$ denote the matrix corresponding to the time series of the delayed variables, i.e.,
\begin{equation*}
    X_{\tau} \coloneqq \begin{bmatrix}
        x_1(t_1-\tau) & x_2(t_1-\tau) & \dots & x_n(t_1-\tau) \\
        x_1(t_2-\tau) & x_2(t_2-\tau) & \dots & x_n(t_2-\tau) \\
        \vdots & \vdots & & \vdots \\
        x_1(t_m-\tau) & x_2(t_m-\tau) & \dots & x_n(t_m-\tau)
    \end{bmatrix}.
\end{equation*}
Note that $X_\tau$ is just a retarded (or even shifted) version of $X$. An augmented library matrix $\Theta(X,X_\tau)\in\mathbb{R}^{m\times q}$ is then introduced by collecting all the basis functions applied to both the current-time and delayed data, e.g.,
\begin{equation*}
\Theta(X,X_{\tau})\coloneqq
\begin{bmatrix}
\;1 & X & X_{\tau} & X^2 & XX_{\tau} & X_{\tau}^{2} & \cdots & \sin(X) & \sin(X_{\tau}) & \cdots \;\;
\end{bmatrix}.\footnote{Note that for a polynomial library of degree $d$, the number of basis functions increases from $p={n+d\choose d}$ in the non-delayed case to $q={(k+1)n+d\choose d}$ in the case of $k$ delays.}
\end{equation*}
The approximation problem now reads
\begin{equation} \label{eq:sindyunknown}
    X'=\Theta(X,X_\tau)\Xi(\tau),
\end{equation}
where the dependence of the sparse solution $\Xi$ from the value of $\tau$ is written explicitly.
Given $\tau$, \eqref{eq:sindyunknown} is solved as \eqref{eq:sindyknown}, i.e., column-wise through
\begin{equation*}
\xi_i(\tau) = \argmin_{\xi \in \mathbb{R}^q}\left(\left\|X_{i}'-\Theta(X,X_\tau)\xi\right\|+\lambda\left\|\xi\right\|_1\right).
\end{equation*}

\section{Bayesian Optimization}
\label{sec:bopt}
BO \cite{Rasmussen06, Shahriari2016, garnett2023} is a class of techniques used to solve the problem
\begin{equation}\label{eq:bo}
\bar{y}=\argmin_{y\in Y}g(y),
\end{equation}
where $g \colon Y \to \mathbb{R}$ with $Y\subseteq\mathbb{R}^{s}$ for some positive integer $s$. It is based on the idea that we can see $g$ as a stochastic process and thus construct a statistical model using its evaluations. BO is applied when no analytical information are known about $g$. Moreover, one assumes that the domain $Y$ has a simple geometry, like a rectangle or a simplex, and typically $s \leq 20$. A classical example of application of BO is the optimization of an expensive-to-evaluate black-box function $g$ without any knowledge of its derivatives nor its properties. In this section we just recall the basic features and properties, in view of combining BO tools with SINDy in Section~\ref{sec:methods} to effectively tackle the problem of identifying unknown delays or other model parameters.

BO is composed of two elements: a probability model $p(g)$ for the objective function $g$ and an acquisition function $a: Y\to \mathbb{R}$ to decide where to sample the objective function. The probability model $p(g)$ captures our current belief of $g$ and is updated into the posterior probability distribution $p(g | D)$ every time a new evaluation $D$ is available. The acquisition function $a$ quantifies for each point in $Y$ the utility of evaluating this point. BO iterates the following three steps (indexed by $i$) until a predefined stopping criterion is fulfilled:
\begin{enumerate}
    \item update the model $p(g)$ on the current data $D_i$;
    \item identify the next point $z_{i+1}\in Y$ on which evaluating $g$ by using the acquisition function:
    \begin{equation*}
        z_{i+1} = \argmin_{y \in Y}a(y);
    \end{equation*}
    \item update the data set as $D_{i+1}\coloneqq D_{i}\cup\{g(z_{i+1})\}$.
\end{enumerate}

The probability model of $g$ is given by a Gaussian process. Every finite collection of random variables composing a Gaussian process has a normal multivariate distribution. These distributions are completely determined by the mean vector and the covariance matrix. To produce the actual model, we first assume to have a set of observations
\begin{equation*}
    D_{\text{obs}}=\{(y_j, u_j)\in Y\times\mathbb{R},\ j = 1,\dots,l \} ,
\end{equation*}
and a number of test points
\begin{equation*}
    D_{\text{test}}=\{ y^*_j\in Y,\ j = 1,\dots,l^*  \}.
\end{equation*}
Then we select a covariance function $k(y,y^*)$ to calculate the covariance matrices. This function tells the statistical relation between different values of the input. One of the most important kernel is the \textit{exponential quadratic} one:
    \begin{equation*}
        k(y,y^*)\coloneqq \sigma^2 \exp \left(-\frac{ \left\Vert y - y^* \right\Vert^2}{2\ell^2}\right),
    \end{equation*}
in which $\sigma$ and $\ell$ are hyperparameters. A list of kernels with a discussion of their properties can be found in \cite{Rasmussen06}.

In the following, we denote with $Y$ the collective input given by the points $y_i$ and with $Y^*$ the collective input given by the points $y_i^*$. We will use $K(Y,Y^*)$ to denote the matrix with entries $k(y_i,y_j^*)$ with $i=1,\dots, l$ and $j=1,\dots,l^*$. We denote with $u=(u_1,\dots,u_l)^T$ the vector given by the observed values and with $u^*=(u_1^*,\dots,u_{l^*}^*)^{T}$ the vector of predicted output for the test set. We can observe that in the context of Gaussian processes what we want is to predict a Gaussian distribution for each test input.

It is possible to prove that $u^*$ has a Gaussian multivariate distribution and it is possible to express explicitly the mean and covariance of this distribution \cite{Rasmussen06}: this is one of the most important advantages of working with Gaussian processes. Assuming that there is no noise in the evaluation of $g$, then $u^*$ is normally distributed with mean vector given by
\begin{equation*}
    \mu = K(Y^*,Y)K(Y,Y)^{-1}u,
\end{equation*}
and covariance matrix given by
\begin{equation*}
    \Sigma = K(Y^*,Y^*) - K(Y^*,Y)K(Y,Y)^{-1}K(Y,Y^*).
\end{equation*}
These equations can be generalized to the case in which the evaluation of $g$ is affected by a random error term \cite{Rasmussen06}. If we focus on a single test point $y^*$, given some observations, we have a prediction, i.e., the maximum a posteriori estimate of the probability distribution of the value of $g$ in $y^*$, and a quantification of the uncertainty, i.e., the standard deviation.

Concerning the acquisition function $a$, the central idea of BO is to trade off the tendency to explore regions in which the model is uncertain and the tendency to exploit the knowledge of the model about good regions in which we have high confidence. An example of acquisition function is the so-called \textit{expected improvement}, where one wants to maximize the expected value over $p(g|D)$:
\begin{equation*}
\mathbb{E}_{p(g|D)}[\max(\overline{u} - g(y), 0)],  
\end{equation*}
where $\overline{u} = \min\{u_1, \dots, u_l\}$ is the current minimum observed value. It has a closed form if $p(g|D)$ is Gaussian \cite{Jones1998}.

The main advantage of using BO to solve optimization problems is that it is a very general technique that uses some natural assumption to model, from a data-driven perspective, a function $g$ and solve the associated optimization problem, providing also a quantification of the uncertainty related to the given solution.

To implement BO in the code we use in Section~\ref{sec:results} we resort to the Python package \texttt{Emukit}~\cite{emukit2023,emukit2019}, a toolkit for emulation and decision making.

\section{Identification of unknown delays and parameters}
\label{sec:methods}
We closed Section~\ref{sec:sindy} with~\eqref{eq:optknown} to solve the optimization problem following the application of SINDy to DDEs with known delays. Knowing the value of the delays in the context of reconstructing a dynamical system from given measurements is rather unrealistic. One could also question whether it is reasonable at all to assume that the system to recover is actually a DDE. This can easily be the case if the available data show complex trends already in the scalar or planar cases (see previous related comments in Section~\ref{sec:dde}) or if the delayed nature of the model is known from available physical knowledge (e.g., a control system or the transmission of a disease). It is thus important to develop a technique capable of identifying systems obtained from DDEs while having no explicit knowledge of the delays. The goal of this section is to devise a strategy to tackle the case of unknown delays.

If the delay $\tau$ is unknown a brute force approach can be used to select among a set $T:=\{\tau^{(1)},\ldots,\tau^{(r)}\}$\footnote{Note in general that the set $T$ can be successively refined around the first optimum found by BO if necessary. Also, by resorting to linear interpolation, there is no restriction on $T$ with respect to the time instants at which the measurement data are available.} of $r$ candidate delay values the value $\bar\tau$ minimizing the reconstruction error of SINDy. This amounts to repeatedly call SINDy on~\eqref{eq:optknown} for every value $\tau\in T$ and then choose
\begin{equation}\label{eq:opttau}
\bar\tau\coloneqq\argmin_{\tau\in T}\left\|X'-\Theta(X, X_{\tau})\Xi(\tau)\right\|.
\end{equation}
This procedure, adopted both in~\cite{sandoz2023sindy} and by the authors in their earlier results presented at \cite{breda2022data}, is highly time-consuming (even though embarrassingly parallelizable). Here instead we propose to use BO to directly solve~\eqref{eq:opttau}, by setting in~\eqref{eq:bo} $y=\tau$, $Y=T$, and $g$ the map defined by
\begin{equation*}
\tau\mapsto\left\|X'-\Theta(X, X_{\tau})\Xi(\tau)\right\| ,
\end{equation*}
through the application of SINDy. In Section~\ref{sec:results} we show that this alternative leads to a  remarkable reduction in the number of calls to SINDy, with a consequent considerable saving in CPU time.

Yet the above saving is not the only advantage of using BO. BO represents an intrinsically multivariate approach, thus the proposed strategy can be applied straightforwardly to the case of multiple unknown delays. It is enough to interpret $\tau$ as a vector $(\tau_{1},\ldots,\tau_{k})^{T}$ of unknown delays, and to optimize $\bar\tau$ in~\eqref{eq:opttau} over $T\coloneqq T_{1}\times\cdots\times T_{k}$ for $T_{j}:=\{\tau_{j}^{(1)},\ldots,\tau_{j}^{(r_{j})}\}$, $j=1,\ldots,k$, with each $r_{j}$ a (possibly different) positive integer. We also note that the proposed procedure is in principle able to determine also the correct (true) number of delays as long as the latter is lower than the value of $k$ considered above.

We conclude this section with two other observations. First, the problem of identifying the delays could be solved by inserting in the augmented
library as much delayed terms $X_{\tau_{j}}$ as possible and leave to a single call of SINDy the selection of  the corresponding model. Although this strategy would require only a single call to SINDy, it gives much less control on the final form of the model and is much more prone to overfitting due to the very large library\footnote{The trivial delay $\tau^{(0)}\coloneqq0$ (as well as other relatively small values) should be excluded from the set $T$ to avoid overfitting, given that $x(t-\tau)\to x(t)$ as $\tau\to0$ thanks to continuity of solutions.}. Second, the multivariate capabilities of the BO strategy can be exploited to identify not only delays, but also any other non-multiplicative parameter. This is the case if one suspects the presence of peculiar nonlinear terms whose functional form depends on parameters, e.g., exponentials $e^{\alpha x}$ or rational terms like $1/(1+x^{\alpha})$ to name a few. This a-priori knowledge is in general driven by physics-informed considerations, knowing that the given measurements come, e.g., from biological surveys or chemical experiments, which often includes typical functional responses such as Michaelis–Menten kinetics or more general Hill functions (see \cite{BLANCHINI2023MICHAELIS} and the references therein). In this sense we will test our approach on the Mackey-Glass equation~\eqref{eq:mg}, trying to simultaneously identify both the delay $\tau$ and the Hill coefficient $\alpha$ of the rational nonlinearity describing the Hill-type response, see~\eqref{eq:hill}.

\section{Numerical results}
\label{sec:results}
This section presents the results of the application of the proposed technique to time series from the four different DDEs described in Section \ref{sec:dde}. The data are obtained by numerical integration of related IVPs through MATLAB's command \texttt{dde23}, with constant initial functions. We present the results in order of complexity, considering first the case of a scalar DDE with a single delay (the logistic DDE of Example~\ref{ex:logistic}), second the case of systems (the delayed SIR model of Example~\ref{ex:sir}), third the case of identifying both delays and parameters (for the Mackey-Glass equation of Example~\ref{ex:mg}), and fourth the case of multiple delays (equation~\eqref{eq:multiple} of Example~\ref{ex:multiple}). From the dynamical point of view, we deal with time series  with either oscillatory or even seemingly chaotic behaviors. The focus is posed both on the correctness of the reconstructed right-hand side and on the reduction of the computational cost given by the BO strategy with respect to the brute force approach described at the beginning of Section~\ref{sec:methods}. Relevant demo codes to reproduce the experimental tests are available at \url{https://github.com/alepec98/SINDy4DDEs}.

\subsection{Delay logistic equation}
\begin{figure}[ht]
\centering
\includegraphics[width=\textwidth]{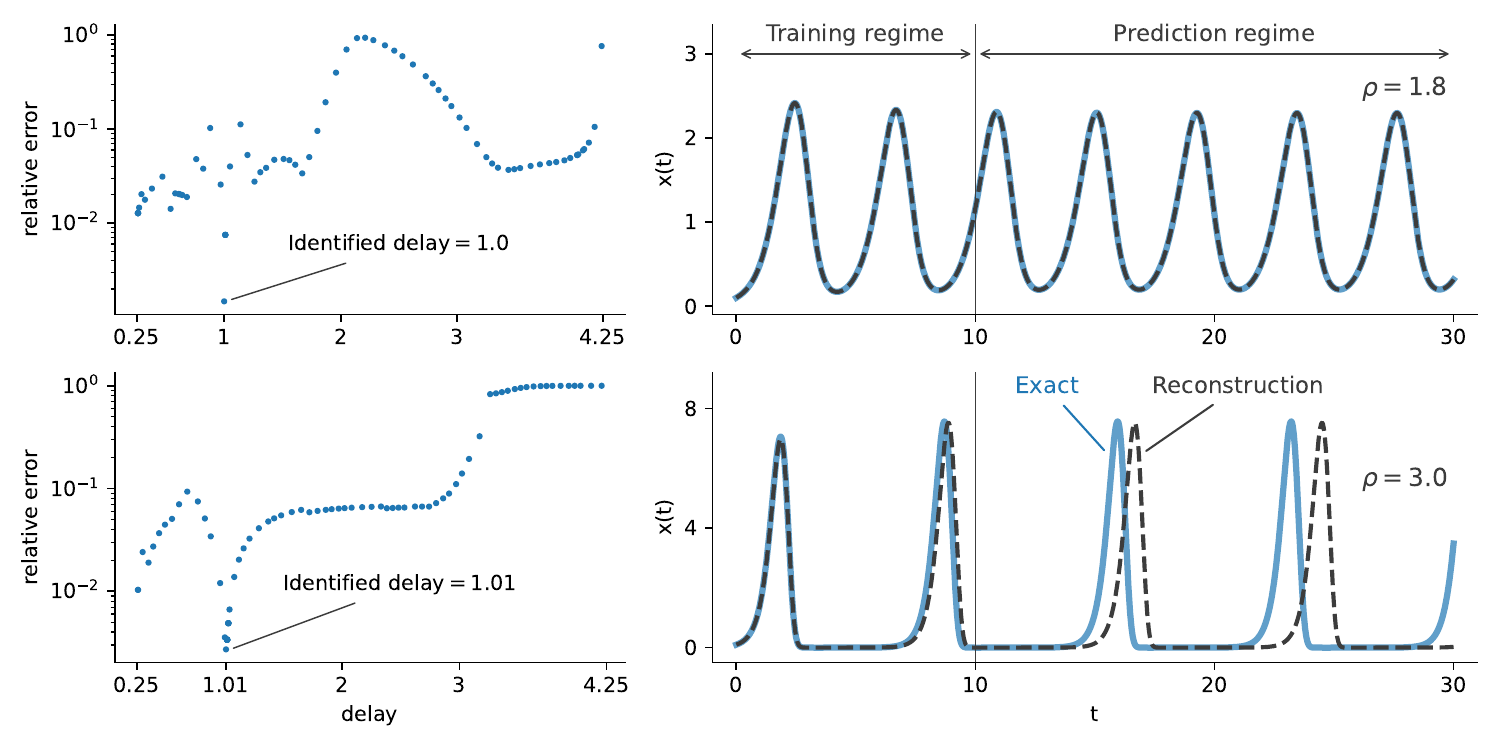}
\caption{Delay logistic equation representative results. The top row shows the test case with $\rho = 1.8$, while in the bottom row we have $\rho = 3.0$. The left column presents the BO evaluations and the identified delays. In the right column we compare the exact solution with the reconstruction corresponding to the identified delay emphasizing the training and the prediction regime.}\label{fig:logistic} 
\end{figure}

We consider as measurement data the time series corresponding to the solutions of the two IVPs for~\eqref{eq:logistic} with parameter values $\rho=1.8$ and $\rho=3.0$, respectively. For both the cases we set the initial function $\varphi=0.1$ and the delay $\tau=1$. The solutions are sampled over the time interval $[0,30]$ with constant stepsize $\Delta t=0.01$. The interval $[0,10]$ is used to train SINDy, while the interval $[10,30]$ is used for validating the recovered DDE via simulation. Assuming the true value of the delay is unknown, we fix the set $T=\{0.25, 0.26,\ldots,4.24\}$ of $r=400$ candidate delays.\footnote{In a practical context, an upper bound for $T$ could be obtained in presence of an oscillatory trend by looking at the frequency of the oscillations.} We choose a polynomial library up to the second degree, giving rise to
\begin{equation*}
\Theta(X,X_{\tau})=
\begin{bmatrix}
1 & X & X_{\tau} & X^2 & XX_{\tau} & X_{\tau}^{2} 
\end{bmatrix}.
\end{equation*}
The results obtained by applying the proposed approach are shown in Figure~\ref{fig:logistic}, first row for $\rho=1.8$, second row for $\rho=3.0$. The left column shows the BO objective function $g=g(\tau)$. The right column shows a comparison between the original time series and the one obtained by solving the identified DDE with the same constant initial function $\varphi$. For $\rho=1.8$ the BO approach correctly identifies a global minimum $\bar\tau=1.0$ and, correspondingly, SINDy extracts from the library the terms $X$ and $XX_{\tau}$, with coefficients $1.8$ and $-1.8$, respectively. In doing so the BO approach requires on average\footnote{For each test case we run the BO with \texttt{Emukit} 10 times, confirming the robustness of the approach and relevant results. The reported values of calls and CPU times are averages over all the runs.} $85.5$ calls to SINDy, corresponding to a $78.6\%$ reduction with respect to a brute force exploration of the whole set $T$. For $\rho = 3.0$ a global minimum is identified again, this time at $\bar\tau=1.01$. Correspondingly, SINDy returns the DDE
\begin{equation*}
x'(t)=2.993x(t)-3.077x(t)x(t-1.01).
\end{equation*}
The small discrepancies (below $3\%$) with respect to the original model justify the deviation between the relevant trajectories, see the bottom-right panel in Figure~\ref{fig:logistic}. The average number of SINDy calls is $96.6$, confirming the $75.9\%$ reduction with respect to the brute force alternative.

\subsection{SIR model with delay}
\begin{figure}[ht]
\centering
\includegraphics[width=\textwidth]{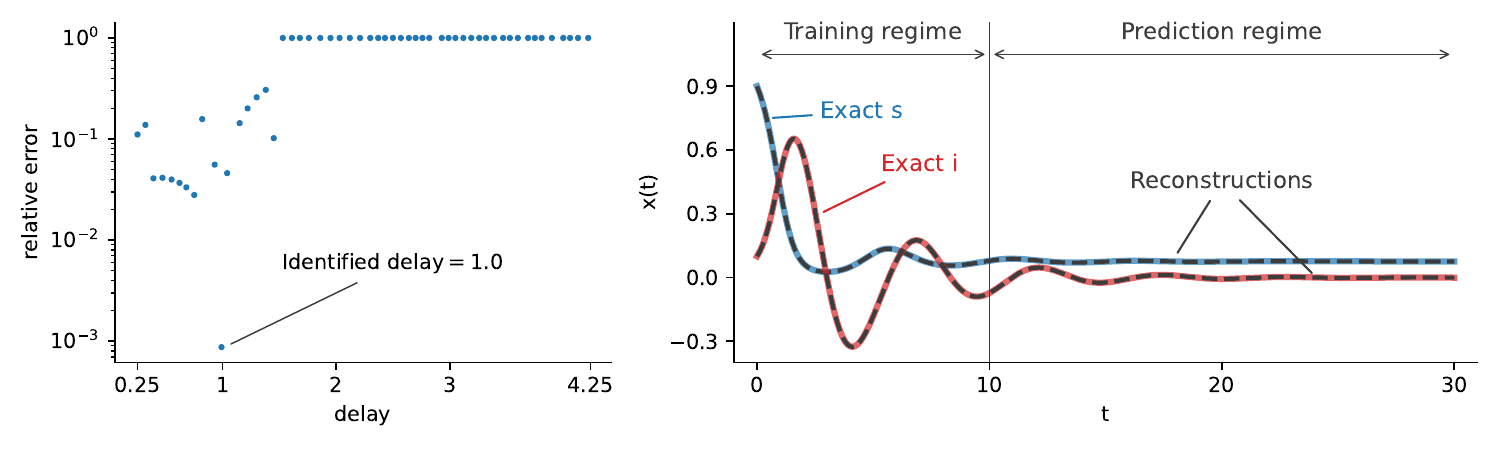}
\caption{SIR model with delay representative results. On the left the BO evaluations and the identified delay. On the right the comparison between the exact solutions and the reconstructions corresponding to the delay $\tau = 1.0$.}\label{fig:sir} 
\end{figure}

As a test case for systems of DDEs we consider~\eqref{eq:sir} in Example~\ref{ex:sir}. The data are taken from simulating the corresponding IVP with constant initial function $\varphi=(0.9,0.1,0)^{T}$, parameters $\beta=3$ and $\mu=1$, and delay $\tau=1$. The solution in $[0,30]$ is sampled at $\Delta t=0.01$ and $T=\{0.25, 0.26,\ldots,4.24\}$. The training interval is $[0,10]$ and the test is done in $[10,30]$. In what follows we consider only susceptible $s$ and infectious $i$ individuals, since recovered ones are obtained as $r=1-s-i$ by assuming constant (normalized) population. Note that including in the SINDy process variables that are linked via algebraic relations could result in different equations --- although equivalent --- from the one we expect\footnote{It is enough to think to this example: reconstruct $\cos (3t)$, with a polynomial library of thirs degree in $\cos (t)$ and $\sin (t)$ over $[0,10]$ considering $\Delta t = 0.1$. SINDy reconstructs the equation $\cos (3t) = -0.667 \cos (t) + 1.667 \cos^3 (t) -2.333 \cos (t)\sin^2(t)$ equivalent but apparently different from the expected $\cos (3t) = 4\cos^3 (t) - 3 \cos (t)$. }. We adopt a polynomial library up to the second degree that combines all samples from $s(t)$, $i(t)$ and $i(t-\tau)$, hence
\begin{equation*}
\Theta(S,I,I_{\tau})=
\begin{bmatrix}
1 & S & I & I_{\tau} & S^2 & SI & SI_{\tau} & I^2 & II_{\tau} & I_{\tau}^2
\end{bmatrix}.
\end{equation*}
The results obtained by applying the proposed approach are shown in Figure~\ref{fig:sir}, with the BO objective function $g=g(\tau)$ (left panel) and a comparison between the original time series and the one obtained by solving the identified system of DDEs with the same constant initial function $\varphi$ (right panel). The BO procedure correctly identifies $\bar\tau=1.0$ and SINDy returns the following DDE:
\begin{equation}
\left\{\setlength\arraycolsep{0.1em}\begin{array}{rcl}
        s'(t)&=&-3.001 s(t)i(t), \\[1mm]
        i'(t)&=&3.001 s(t)i(t)- 1.0 i(t-1.0).
\end{array}\right.
\end{equation}
On average the SINDy algorithm is called $69.3$ times, thus achieving $82.7\%$ less evaluations of the error function.

\subsection{Mackey-Glass equation}
\begin{figure}[ht]
\centering
\includegraphics[width=\textwidth]
{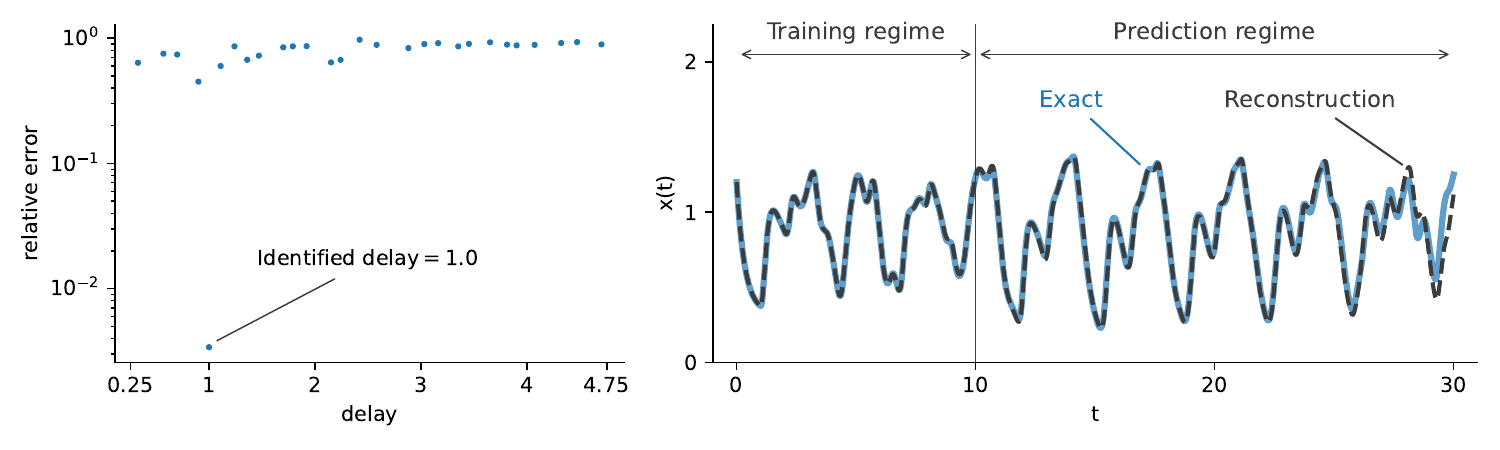}
\caption{Mackey-Glass equation results for unknown delay $\tau$. On the left the BO evaluations and the identified delay.
On the right the comparison between the exact solution and the reconstruction corresponding to
the delay $\tau = 1.0$.}\label{fig:mg1} 
\end{figure}

Data for the Mackey-Glass equation~\eqref{eq:mg} of Example~\ref{ex:mg} are obtained by simulating the corresponding IVP with constant initial function $\varphi=0.1$, parameters $\beta=6$, $\gamma=3$, $\alpha=10$, and delay $\tau=1$. The solution in $[0,30]$ is sampled with $\Delta t = 0.01$. The training is done using data from the interval $[0,17]$ while the interval $[17,30]$ is kept as test. As observed in Section \ref{sec:dde}, the right-hand side of the Mackey-Glass equation is not restricted to polynomial terms, due to the presence of the Hill response~\eqref{eq:hill}. This class of functional feedback is common in biological modeling, thus we feel allowed to include it in the SINDy library in a physics-informed fashion. Nevertheless, in a first test we consider as unknown only the delay $\tau$, using the candidate set $T=\{0.25, 0.26,\ldots,4.74\}$, hence $r=450$. The SINDy library is based on second-degree polynomials and reads
\begin{equation}\label{eq:libmg}
\Theta(X,X_{\tau},H_{\tau;\alpha})=
\begin{bmatrix}
1 & X & X_{\tau} & H_{\tau;\alpha} & X^2 & XX_{\tau} & XH_{\tau;\alpha} & X_{\tau}^2 & X_{\tau}H_{\tau;\alpha} & H_{\tau;\alpha}^2
\end{bmatrix} ,
\end{equation}
with $\alpha=10$ assumed to be known and $H_{\tau;\alpha}$ the sampling of \eqref{eq:hill}. The results of the proposed technique are visible in Figure~\ref{fig:mg1}. BO correctly identifies a global minimum $\bar\tau=1.0$ of the objective function. Correspondingly, SINDy returns the following DDE:
\begin{equation}\label{eq:mgsindy}
    x'(t) =  5.997\frac{ x(t-1.0) }{1+(x(t-1.0))^{10}}-2.998 x(t).
\end{equation}
The average number of BO iterations is $50.9$, providing a reduction of $88.7\%$ of total SINDy calls with respect to the brute force approach.

\begin{figure}
\centering
\includegraphics[width=\textwidth]{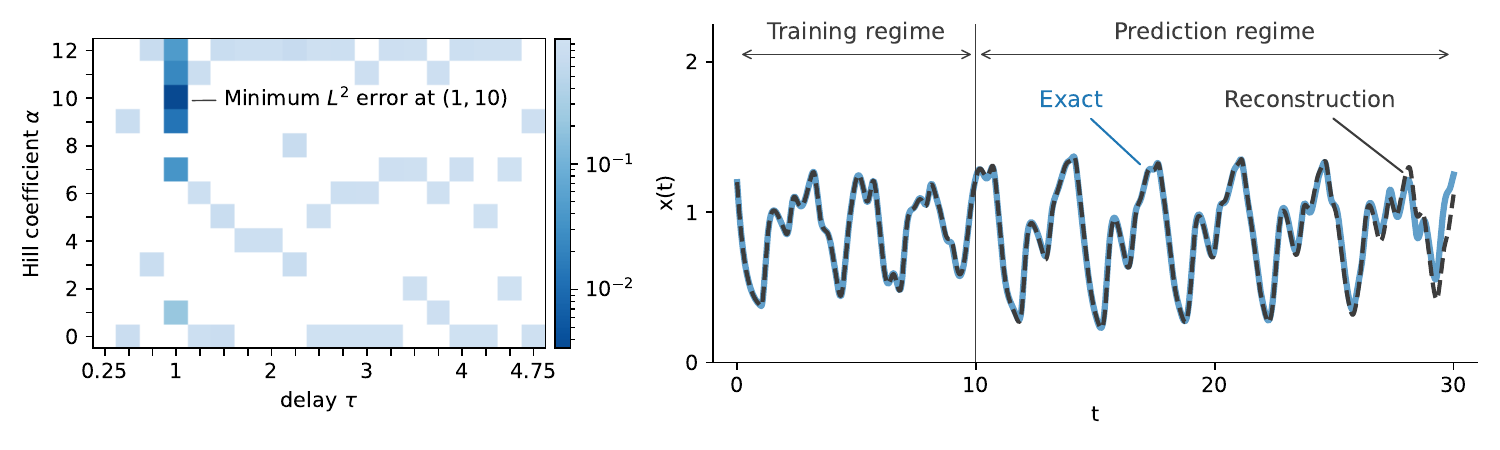}
\caption{Mackey-Glass equation representative results for unknown delay and Hill coefficient. On the left the BO evaluations and the identified couple of delay and Hill coefficient.
On the right the comparison between the exact solution and the reconstruction corresponding to
$\tau = 1.0$ and $\alpha = 10$.}\label{fig:mg2} 
\end{figure}

In a second test we consider as unknown both the delay $\tau$ and the Hill coefficient $\alpha$, using two different candidate sets, namely $T=\{0.25,0.50,\ldots,4.75\}\times\{0,1,2,\dots,12\}$ (made of $19\times13=247$ values) and $T=\{0.25,0.30,\ldots,4.75\}\times\{0,1,2,\dots,12\}$ (made of $91\times13=1183$ values). The SINDy library is always~\eqref{eq:libmg}, this time with unknown $\alpha$. Note that now $y$ in~\eqref{eq:bo} is interpreted as the couple $(\tau,\alpha)$. The training data are the same as for the first case. An illustration of the results can be seen in Figure~\ref{fig:mg2}. In both cases, BO is able to identify the correct minimum $(\bar\tau,\bar\alpha)=(1.0,10)$ and SINDy returns again Equation~\eqref{eq:mgsindy}. The average number of BO iterations is $132.5$ (corresponding to a reduction of $46.4\%$ of SINDy calls) for the first candidate set and $224.1$ (corresponding to a reduction of $81.1\%$ of SINDy calls) for the second candidate set. Figure~\ref{fig:mg2} shows, for the first candidate set, the BO evaluations with the corresponding reconstruction errors (left panel) and a comparison between training and reconstructed trajectories (right panel).

\subsection{Multiple delays}
\begin{figure}[ht]
\centering
\includegraphics[width=\textwidth, trim=50 0 0 0, clip]{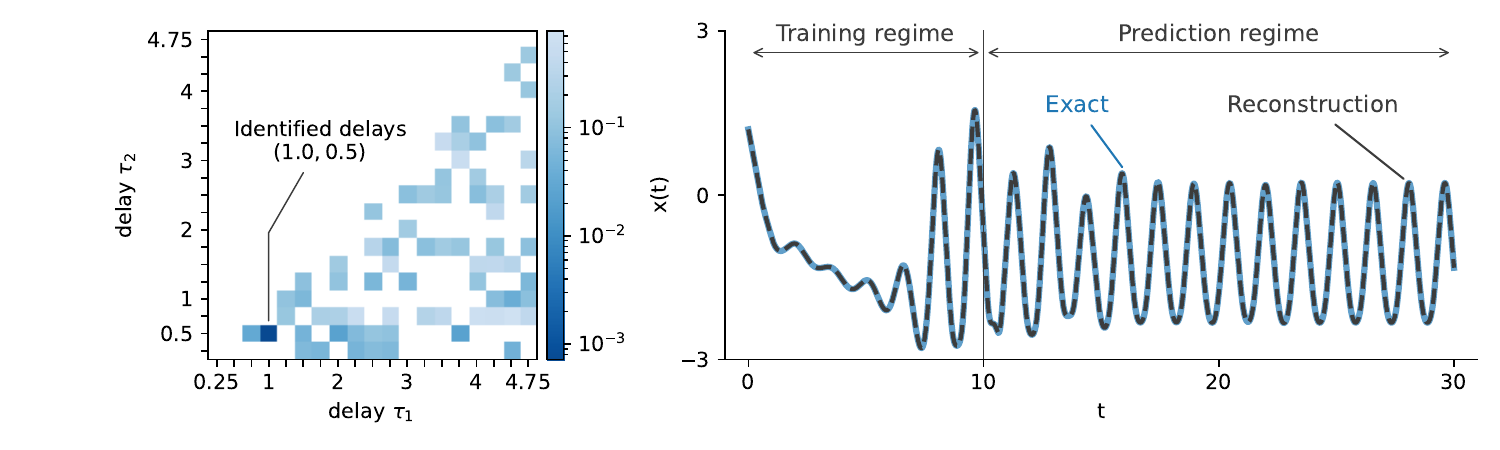}
\caption{Multiple delays representative results. On the left the BO evaluations and the identified delays.
On the right the comparison between the exact solution and the reconstruction corresponding to
$\tau_1 = 1.0$ and $\tau_2 = 0.5$.}\label{fig:multiple}
\end{figure}

For the multiple delays test case, data are obtained by simulating the corresponding IVP for~\eqref{eq:multiple} in Example~\ref{ex:multiple} with initial function $\varphi=1.2$, parameters $a=-1.0$, $b=-0.5$, and delays $\tau_{1}=1.0$ and $\tau_{2}=0.5$. Uniform samples on $[0,30]$ are extracted with $\Delta t=0.01$ and used for training in $[0,10]$ and test in $[10,30]$. Based on a third-degree polynomial library, SINDy returns the following DDE:
\begin{equation}\label{eq:multiplesindy}
x'(t)=-1.0 x^2 (t-1.0) - 0.5 x^3(t-0.5) ,
\end{equation}
for the case of known delays. When the delays are assumed to be unknown, the BO strategy is applied for three candidate sets, namely $T=\{0.25,0.50,\ldots,4.75\}\times\{0.25,0.50,\ldots,4.75\}$ (made of $19\times19=361$ values), $T=\{0.20,0.30,\ldots,4.20\}\times\{0.20,0.30,\ldots,4.20\}$ (made of $41\times41=1\,681$ values) and $T=\{0.25,0.30,\ldots,3.25\}\times\{0.25,0.30,\ldots,3.25\}$ (made of $61\times61=3\,721$ values). We enforced the constraint $\tau_{1}>\tau_{2}$ without loss of generality, so that the total number of candidate delays is roughly halved ($171$ for the first case, $820$ in the second case, and $1\,830$ in the third case). Figure~\ref{fig:multiple} shows the numerical results. In all the cases, BO is able to identify the correct minimum $(\bar\tau_{1},\bar\tau_{2})=(1.00,0.50)$ and SINDy returns always~\eqref{eq:multiplesindy}. The average number of BO iterations is $101.7$ (corresponding to a reduction of $40.5\%$ of SINDy calls) for the first set, $510.4$ (corresponding to a reduction of $37.8\%$ of SINDy calls) for the second candidate set and $384.0$ (corresponding to a reduction of $79.0\%$ of SINDy calls) for the third candidate set. Figure~\ref{fig:multiple} shows, for the first candidate set, a color map of the objective function (left panel) and a comparison between training and reconstructed trajectories (right panel). Note that the lower part of the outer panel displays the maximum error due to the forced constraint on the delays.

\subsection{Discussion}
\label{sec:discussions}
The proposed methodology achieves good results in all the considered test cases. We were able to correctly identify the governing equations from data with minimal errors, even in the case of unknown delays and non-multiplicative parameters through the application of BO in place of a brute force approach. Moreover, BO remarkably reduces the computational burden with an average reduction of the number of evaluations of the error function equal to $67.9\%$.

The case $\rho = 1.8$ of the logistic equation~\eqref{eq:logistic} is reconstructed perfectly, with exact correspondence between the original coefficients and delays and the identified ones. In the case $\rho = 3.0$, the delays and the coefficients are determined with an error smaller than $3\%$, providing a good reconstruction considering the fact that we are in a data-driven setting. The small error in the delay is the cause of the discrepancy we see in the plot in Figure~\ref{fig:logistic} (bottom-right panel). In particular, the frequency of the peaks depends on the limit cycle, whose structure can be very different for two close delays.

In the analysis of the SIR model with delay~\eqref{eq:sir} we saw a good reconstruction of the original equations for $s$ and $i$. We emphasize how the coefficient $\beta$, which is present in both the equations, is perfectly reconstructed. This behaviour can be also seen in the reconstruction of the logistic equation with $\rho = 1.8$. In this case, the coefficient $\rho$ appears twice if the right-hand side is written explicitly. However, SINDy identified the same values for the two different occurrences.

The study of the Mackey-Glass equation~\eqref{eq:mg} gives an example of how our methodology can be applied to DDEs with complex structures and chaotic behavior. The delay is identified perfectly and the coefficients present only a negligible error. Even if for training our model we are considering a limited amount of data, the reconstruction remains close to the original solution for a significantly longer time interval. Moreover, our analysis shows that BO can be used, in combination with SINDy, to select the member of a given family of functions depending on a small number of parameters without testing all the members of the family.

Finally, the same good results are obtained in the case of multiple delays as shown for equation~\eqref{eq:multiple}.

\section{Conclusions and future perspectives}
\label{sec:conclusion}
In this work, we proposed a numerical method to identify the governing equations in the case of DDEs and systems of DDEs with unknown discrete delays, and non-multiplicative parameters. Our method is based on a combination of the SINDy algorithm to identify the coefficients of the governing equation and the use of BO to select the correct values for the unknown delays and/or parameters. We tested our method by considering equations and systems with different long-term behavior. The identified coefficients, delays, and parameters were close to the original ones in all the cases considered. Adopting BO to identify the unknown values through the optimization of the error function enables an improvement in performance: BO requires usually $67.9\%$ less evaluations of the error function with respect to the brute force approach of testing the whole set of candidates.

Future works will investigate different extensions of the proposed method. 
A possible research direction could be to study how the choice of the error function affects the system identification, e.g., by considering the error on the state variable instead of on its derivative. This method would help in cases where we have a sensitive dependence on the coefficients and the delays: this is the case of the logistic equation for $\rho = 3.0$ and of equations with chaotic behavior.

Extending the proposed techniques to the case of DDEs with distributed delays is another possibility. In this case one would have to identify not a single delay, but a set of continuous delays together with the kernel function in the integral. An extension to this kind of equations would open different possibilities of applications to the study of structured populations from the data-driven point of view.

\section*{Acknowledgements}
This work was partially supported by European
Union Funding for Research and Innovation --- Horizon 2020 Program --- in the
framework of European Research Council Executive Agency: H2020 ERC CoG 2015
AROMA-CFD project 681447 ``Advanced Reduced Order Methods with Applications in
Computational Fluid Dynamics'' P.I. Professor Gianluigi Rozza and by the GNCS 2023 project ``Sistemi dinamici e modelli di evoluzione: tecniche funzionali, analisi qualitativa e metodi numerici'' (CUP: E53C22001930001).

\bibliographystyle{abbrvurl}

\end{document}